\documentclass[letterpaper, 10 pt, conference]{ieeeconf}  

\IEEEoverridecommandlockouts                              

\pdfoutput=1

\usepackage{algorithm}
\usepackage{algpseudocode}
\usepackage{amsmath,amsfonts,amssymb}
\usepackage{graphics} 
\usepackage{epsfig,color} 
\usepackage{times} 
\usepackage{amssymb}  
\usepackage{amsthm}
\usepackage{graphicx}
\graphicspath{ {./images/} }
\usepackage{graphicx}
\usepackage{subfigure}
\usepackage{booktabs} 
\usepackage{caption}

\usepackage{subcaption}
\usepackage{multirow}
\usepackage{bbm}

\theoremstyle{definition}
\newtheorem{definition}{Definition}[section]

\newcommand{\Bcal}{{\cal B}}

\newcommand{\Fcal}{{\cal F}}

\newcommand{\Hcal}{{\cal H}}

\newcommand{\Ocal}{{\cal O}}
\newcommand{\Pcal}{{\cal P}}

\newcommand{\Tcal}{{\cal T}}

\newcommand{\argmin}{\mathop{\rm argmin}}

\newtheorem{thm}{Theorem}

\newtheorem{assump}{Assumption}
\newtheorem{remark}{Remark}

\usepackage[bookmarks=true,pageanchor,colorlinks,linkcolor=red,anchorcolor=blue, citecolor=blue,urlcolor=blue,hyperfootnotes=false]{hyperref}
\title{\LARGE \bf
Resilient Federated Learning under Byzantine Attack in Distributed Nonconvex Optimization with $2$-$f$ Redundancy

}

\author{ Amit Dutta \quad Thinh T. Doan\quad Jeffrey H. Reed 
\thanks{The authors are with the Bradley Department of Electrical and Computer Engineering, Virginia Tech, Blacksburg, VA. Email: {\tt\small \{amitdutta, thinhdoan, reedjh\}@vt.edu}}%
}

\begin{document}

\maketitle
\thispagestyle{empty}
\pagestyle{empty}

\begin{abstract}
We study the problem of Byzantine fault tolerance in a distributed optimization setting, where there is a group of $N$ agents communicating with a trusted centralized coordinator. Among these agents, there is a subset of $f$ agents that may not follow a prescribed algorithm and may share arbitrarily incorrect information with the coordinator. The goal is to find the optimizer of the aggregate cost functions of the honest agents. We will be interested in studying the local gradient descent method, also known as federated learning, to solve this problem. However, this method often returns an approximate value of the underlying optimal solution in the Byzantine setting. Recent work showed that by incorporating the so-called comparative elimination (CE) filter at the coordinator, one can provably mitigate the detrimental impact of Byzantine agents and precisely compute the true optimizer in the convex setting. The focus of the present work is to provide theoretical results to show the convergence of local gradient methods with the CE filter in a nonconvex setting. We will also provide a number of numerical simulations to support our theoretical results.

\end{abstract}


\section{Introduction}

The constant expansion of large networks has resulted in an exponential growth of data, creating a pressing need for higher computational power and storage requirements. To address this demand, numerous distributed algorithms have been developed, where computation and data are distributed over networks. Federated learning has emerged as a popular distributed framework that facilitates collaborative training of a shared model by multiple devices \cite{mcmahan2017communication, zhao2021fedpage, woodworth2020local}. This technique is particularly relevant in the context of wireless communication networks \cite{niknam2020federated}, \cite{yang2020energy}, where there is a growing need for efficient and scalable machine learning solutions due to the rapidly increasing number of wireless devices. In this context, a common problem often reduces to optimize an aggregate objective function that is composed of $N$  functions distributed at $N$ different  agents (e.g., mobile devices). 
By using federated learning framework, in applications like  federated spectrum sensing over a wireless sensor network, the updates of any optimization algorithm can be implemented locally at the agents without requiring data being transmitted to the centralized coordinator (e.g., a server/network operator). By keeping data locally, federated learning not only reduces the amount of communications between agents and the server but also provides some level of privacy.  

One of the main challenges in federated learning is the vulnerability of the system to malicious attacks where some agents in the network may fail or whose updates can be manipulated by an external entity. Such malicious (Byzantine) agents will have detrimental impacts to the performance of other agents, and if not addressed, it can lead to catastrophic failures of entire network \cite{shi2021challenges}.  

In this paper, we study the performance of federated learning, in particular, the celebrated distributed local stochastic gradient descent (SGD), when a (small) number of agents in the network is malicious. We will focus on Byzantine malicious attacks, where Byzantine agents can observe the entire network and send any information to the centralized coordinator. Under the presence of Byzantine agents, it is impossible for nonfaulty agents to find the optimizer of the aggregate of functions at every agent (including Byzantine agents) as Byzantine agents can send a random number irrelevant to its function. Thus, we consider another meaningful objective in this setting, where the goal is to solve the optimization problem only involving the honest agents. In particular, we consider the setting where there are up to $f$ faulty Byzantine agents with unknown identities. Our goal is address the following exact fault-tolerance problem.

\textbf{Exact fault-tolerance problem}: 
Let $\mathcal{H}$ be the set of honest agents with $|\mathcal{H}| \geq N-f$. A distributed optimization algorithm is said to have exact fault-tolerance if it allows all the non-faulty agents to compute
\begin{align}\label{eq:exact_tolerence}
    x_{\mathcal{H}}^{*} \in \argmin_{x\in \mathbb{R}^{d}}\sum_{i \in \mathcal{H}}q^{i}(x).
\end{align}
We will study this exact fault-tolerance problem under the following $2f$-\textit{redundancy} condition, which is necessary and sufficient for solving problem \eqref{eq:exact_tolerence}  \cite{gupta2020fault}. 
\begin{definition}[$2f$-\textbf{redundancy}]
The set of non-faulty agents $\mathcal{H}$, with $|\mathcal{H}| \geq N-f$, is said to have $2f$-redundancy if for any subset $\mathcal{S} \subset \mathcal{H}$ with $|\mathcal{S}| \geq N -2f$,
\begin{align}
    \argmin_{x\in \mathbb{R}^{d}}\sum_{i \in \mathcal{S}}q^{i}(x) = \argmin_{x\in \mathbb{R}^{d}}\sum_{i \in \mathcal{H}}q^{i}(x). 
\end{align}
\end{definition}
The definition above states that solving the optimization problem over the honest agents is equivalent to solving it over a minimum of $N-2f$ honest agents. We will study the performance of distributed local SGD under the 2f-redundancy condition. Motivated by the recent work in \cite{gupta2023byzantine}, we will consider the so-called comparative elimination (CE) filter in the distributed local SGD to mitigate the detrimental impact of Byzantine agents. Our focus is to provide theoretical results to show the convergence of distributed local SGD with CE filter in a nonconvex setting.


\noindent\textbf{The main contribution} of this paper is to study the performance of federated local SGD with the CE filter for solving distributed nonconvex optimization problems under Byzantine attacks with the 2f-redundancy condition.  We show that this method solves the exact fault-tolerance problem at a linear rate when the objective satisfies Polyak-\L ojasiewicz (PL) condition. We will also provide a number of numerical simulations to illustrate our theoretical results.

\subsection{Related work}

Existing literature offers various Byzantine fault-tolerant aggregation schemes, such as \textit{multi-KRUM} \cite{blanchard2017machine}, \textit{CWMT} \cite{su2019finite}, \textit{GMoM} \cite{chen2017distributed}, \textit{MDA} \cite{guerraoui2018hidden}, and \textit{Byzantine-RSA} \cite{li2019rsa} filters. However, these schemes don't guarantee exact fault-tolerance without additional assumptions. \cite{gupta2023byzantine}  showed the possibility of achieving exact fault tolerance in a deterministic setting and approximate fault tolerance in a stochastic setting with $2f$ redundancy. Recently, \cite{farhadkhani2022byzantine} proposed RESAM, a unified Byzantine fault-tolerant framework based on previous methods, demonstrating finite-time convergence with additional assumptions. Notably, their results apply to non-convex objectives but exclude the CE aggregation scheme.

Our work in this paper extends work by \cite{gupta2023byzantine} to the non-convex setting, where the global objective function satisfies the PL condition. This broadens the analytical framework beyond strongly convex scenarios. Additionally, relevant work by \cite{cao2019distributed, cao2020distributed} addresses approximate fault tolerance with more relaxed conditions on Byzantine agents.

\section{Federated local SGD with Byzantine agents}\label{sec:CE_Filter}

The proposed federated local SGD with CE filter is formally presented in Algorithm \ref{alg:cap}. 
Each user maintains a local variable $x^{i}_{k,t}$, estimating the local optimal solution at the $k^{\text{th}}$ global iteration and $t^{\text{th}}$ local iteration. The server maintains the global optimal solution estimate $\Bar{x}_{k}$. The server initializes each user by transmitting an initial global model, $\Bar{x}_{0}$. Each user $i$ then runs $\mathcal{T}$ local SGD steps as shown in line $6$, using a step-size $\alpha_{k}$ and a sample $g^{i}(x^{i}_{k,t})$ of its local gradient to update $x^{i}_{k,t}$. In the stochastic setting $g^{i}(x^{i}_{k,t}) = \nabla Q^{i}(x^{i}_{k,t};\Delta^{i}_{k,t})$, where the gradients are i.i.d. sampled. After the local SGD steps, users transmit their estimates to the server. Using the CE filter the server removes all the values that are suspiciously large as compared to its value.  Lines $10$ and $11$ show the CE filter where the server sorts the distance of the user estimates from its average and eliminates $f$ largest distances. Finally, the server then computes a new average based on the estimates of remaining $N-f$ agents as shown in line $12$. 
\begin{algorithm}[!t]
\caption{Federated Local SGD with CE Filter}\label{alg:cap}
\begin{algorithmic}[1]
\State \textbf{Initialize}: The server initializes the model with $\Bar{x}_{0}\in \mathbb{R}^{d}$. Each agent initializes with step-sizes $\alpha_{k}$ and chooses $\mathcal{T}$.
\For{$k=0,1,..$}
    \State All clients $i=1,2,..,N$ in parallel do
    \State Receive $\Bar{x}_{k}$ from the server and set $x^{i}_{k,0} = \Bar{x}_{k}$
     \For{$t=0,...,\mathcal{T}-1$}
     \State $x^{i}_{k,t+1} = x^{i}_{k,t} - \alpha_{k} g^{i}(x^{i}_{k,t}).$
    \EndFor
    \State Users send $x^{i}_{k,T}$ to the server
    \State Server sorts these values as
    \State $\|\Bar{x}_{k} - x^{i_1}_{k,1}\| \leq \|\Bar{x}_{k} - x^{i_2}_{k,1}\| \leq \cdots \leq \|\Bar{x}_{k} - x^{i_N}_{k,1}\|,$
    \State\quad\quad $\Bar{x}_{k+1} = \frac{1}{|\mathcal{F}_{k}|}\sum_{i\in \mathcal{F}_{k}}^{}x^{i}_{k,\mathcal{T}}.$
\EndFor
\end{algorithmic}

\end{algorithm}


\section{Main results}\label{sec:Main_results}

In this section we will provide our theoretical findings on the convergence properties of the local SGD with CE filter. For this we first define the average cost of the non-faulty agents defined as 

\begin{align}
    q^{\mathcal{H}}(\Bar{x}_{k}) = \frac{1}{|\mathcal{H}|}\sum_{i \in \mathcal{H}}q^{i}(x),
\end{align}
Our main results are studied based on some fairly standard assumptions for non-convex optimization as stated below.
\begin{assump}\label{eq:Assump_1}
(Lipchitz smoothness). For each $i \in \mathcal{H}$, $q^{i}$ has L-Lipschitz continuous gradient
\begin{align}\label{eq:Lipschitz}
    q^{i}(y) - q^{i}(x) \leq \nabla q^{i}(x)^{T}(y-x) + \frac{L}{2}\|y-x\|^{2}.
\end{align}
\end{assump}
\begin{assump}\label{eq:Assump_2} The function $q^{\mathcal{H}}$ satisfies PL condition and quadratic growth with some $\mu\geq0$
\begin{align}\label{eq:PL condition}
    \|\nabla q^{\mathcal{H}}(\Bar{x}_{k})\|^{2} \geq 2\mu (q^{\mathcal{H}}(\Bar{x}_{k}) - q^{\mathcal{H}}(x^{*}_{\mathcal{H}})) \geq \mu^{2}\|\Bar{x}_{k} - x^{*}_{\mathcal{H}}\|^{2}.
\end{align}

\end{assump}
\begin{assump}\label{eq:Assump_3}
    The The random variables $\Delta^{i}_{k}$, for all i and k, are i.i.d., and there exists a positive constant $\sigma$ such that
    \begin{align}
        &\mathbb{E}[\nabla Q^{i}(x,\Delta^{i}_{k,t})|\mathcal{P}_{k,t}] = \nabla q^{i}(x), \quad \forall x \in \mathbb{R}^{d},\notag\\
        &\mathbb{E}[\|\nabla Q^{i}(x,\Delta^{i}_{k,t})-\nabla q^{i}(x)\|^{2}|\mathcal{P}_{k,t}] \leq \sigma^{2}, \quad \forall x \in \mathbb{R}^{d}. \label{eq:sigma}
    \end{align}
\end{assump}
In Assumption \ref{eq:Assump_3}, $\mathcal{P}_{k,t}$ is defined as a filtration containing all the history generated by Algorithm \ref{alg:cap} up to time $k+t$. 
\begin{align*}
    \mathcal{P}_{k,t} = \cup_{i \in \Hcal}\{\bar{x}_{0}, \cdots,\bar{x}_{k}, x^{i}_{k,1}, \cdots, x^{i}_{k,t}\}
\end{align*}
Further we have $|\Bcal_{k}|+|\Hcal_{k}| = |\Fcal_{k}|=|\Hcal|$ for any $k \geq 0$. 


Next, we present our main theoretical result of this paper, where we study the convergence rate of Algorithm \ref{alg:cap} in solving problem \eqref{eq:exact_tolerence}. For an ease of exposition, we present the proof of our result in Section \ref{appendix}. 

From step $6$ in Algorithm \ref{alg:cap}, for $i \in \mathcal{H}$ the local update is equivalent to
\begin{align}\label{eq:Algo_T1}
    x^{i}_{k,t+1} = \Bar{x}_{k} - \alpha_{k}\sum_{l=0}^{t}\nabla Q^{i}(x^{i}_{k,l}, \Delta^{i}_{k,l}) .
\end{align}
Our result for the stochastic setting is presented below.

\begin{thm}\label{thm:Stochastic_Theorem2}
Let $\{\Bar{x}_{k}\}$ be generated by Algorithm \ref{alg:cap}. Let $\alpha_{k}$ be chosen as
\begin{align}\label{eq:assumption_2_Theorem2}
    \alpha_{k} = \alpha \leq \frac{\mu}{72L^{2}\mathcal{T}}\cdot
\end{align}
Then, if the following condition holds 
\begin{align}\label{eq:assumption_1_Theorem2}
    \frac{f}{N-f} \leq \frac{\mu}{3L},
\end{align}
then we have
\begin{align}
    &\mathbb{E}[q^{\mathcal{H}}(\Bar{x}_{k+1}) - q^{\mathcal{H}}(x^{*}_{\mathcal{H}})]  \notag \\
    &\leq \Big(1-\frac{\alpha_{k}\mu\mathcal{T}}{36}\Big)^{k+1}\mathbb{E}[q^{\mathcal{H}}(\bar{x}_{0}) - q^{\mathcal{H}}(x^{*}_{\mathcal{H}})] \notag\\
    &\quad +\frac{180L\Tcal\alpha_{k}\sigma^{2}}{\mu}  + \frac{72\Tcal\sigma^{2}f}{\mu|\Hcal|}.\label{thm:Stochastic_Theorem2:ineq}
\end{align}
\end{thm}
\begin{remark}
In Theorem \ref{thm:Stochastic_Theorem2} due to the constant step size, the optimality error converges linearly only to a ball centered at the origin. The size of the ball is determined by two factors. The first factor is dependent on the step size $\alpha$, which is commonly observed in the convergence of local gradient descent with non-faulty agents. The second factor is influenced by the level of gradient noise, denoted by $\sigma$. This noise is a result of both the Byzantine agents and the stochastic gradient samples.

It is worth noting that our comparative filter is specifically designed to eliminate potentially erroneous values sent by the Byzantine agents. However, it is unable to address the issue of variance in their stochastic samples. One possible solution to this problem is to have each agent sample a mini-batch of size $m$, thereby replacing $\sigma^2$ in \eqref{thm:Stochastic_Theorem2:ineq} with $\sigma^2/m$. By increasing the size of $m$, the optimality error can be made arbitrarily close to zero. Furthermore, when $\alpha_k\sim 1/k$, the convergence rate is $\Ocal(1/k)$. 

Finally, if we can have access to the exact values of the gradients $\nabla q^{i}$, then $\bar{x}_{k}$ converges exactly to $x^{\star}$ exponentially. 

\end{remark}



\section{Simulations}\label{sec:simulations}
\begin{figure*}[]\label{sin_1}
    \centering
     \subfigure{
        \includegraphics[width=0.45\linewidth, height=0.25\textwidth]{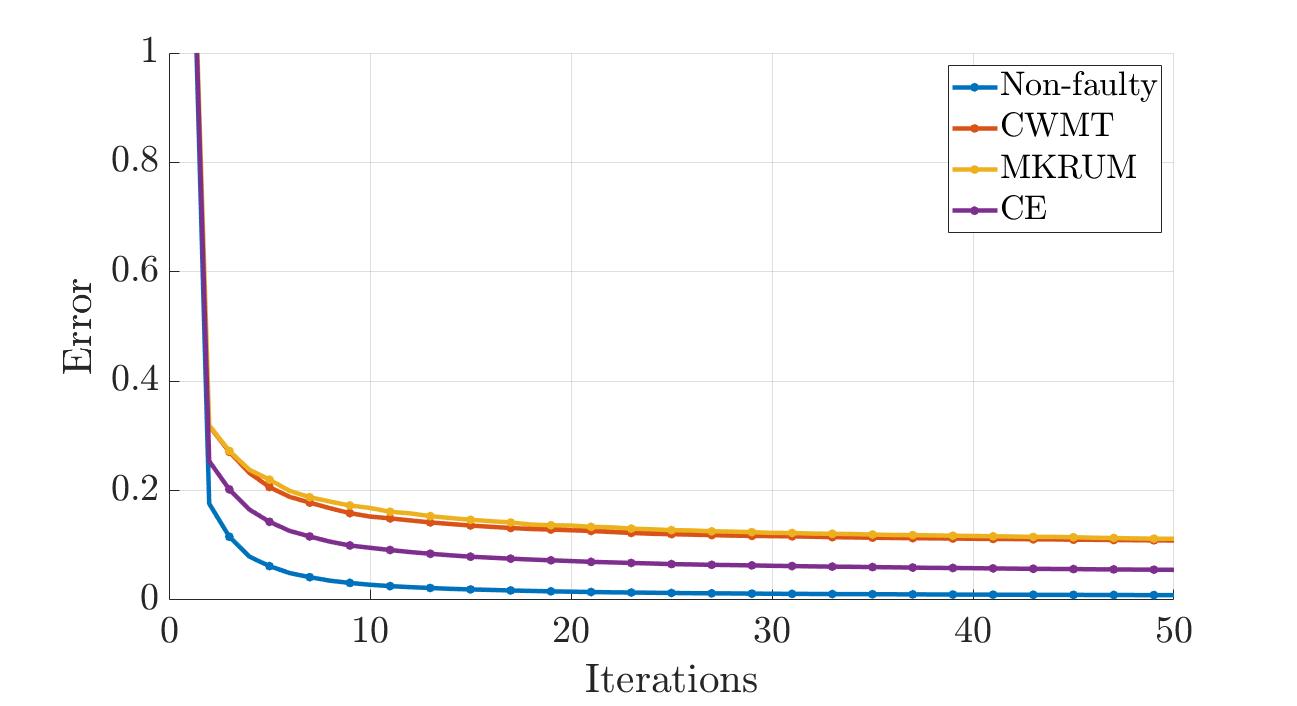}
        \label{fig:gull}}
     \subfigure{
        \includegraphics[width=0.45\linewidth, height=0.25\textwidth]{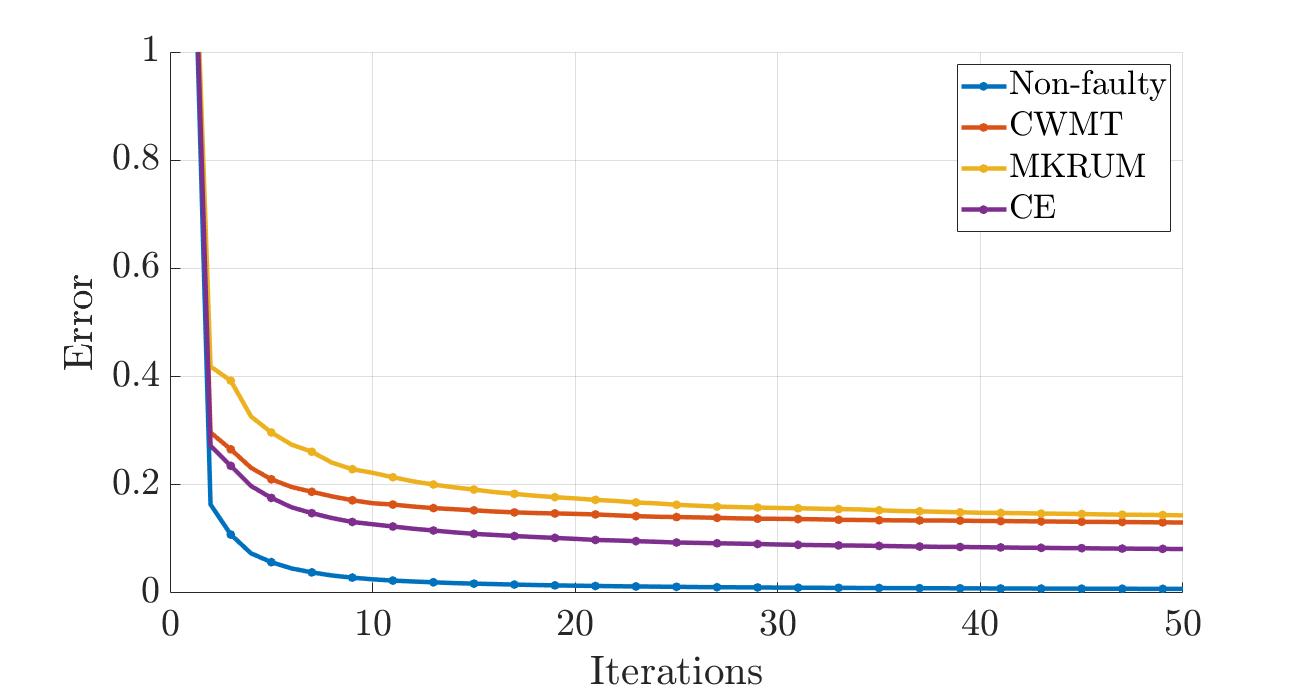}
        \label{fig:tiger}}
        \subfigure{
   \includegraphics[width=0.45\linewidth, height=0.25\textwidth]{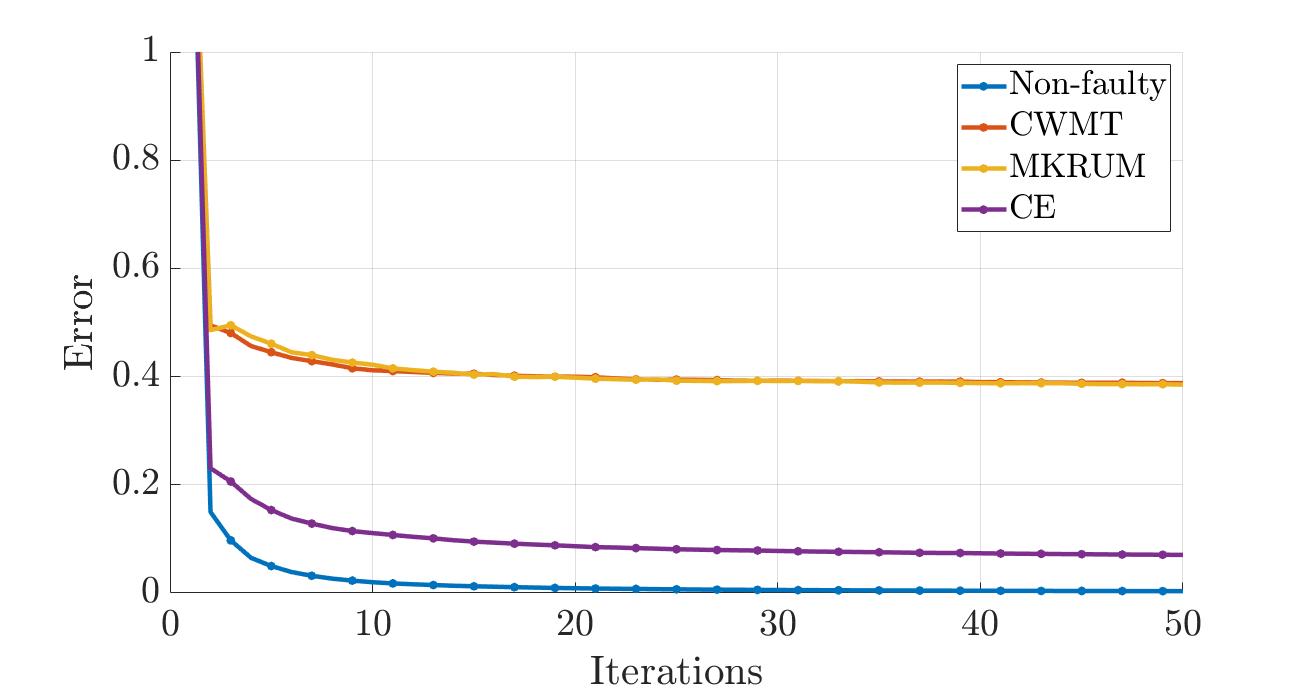}
        \label{fig:tiger}}
        \subfigure{
   \includegraphics[width=0.45\linewidth, height=0.25\textwidth]{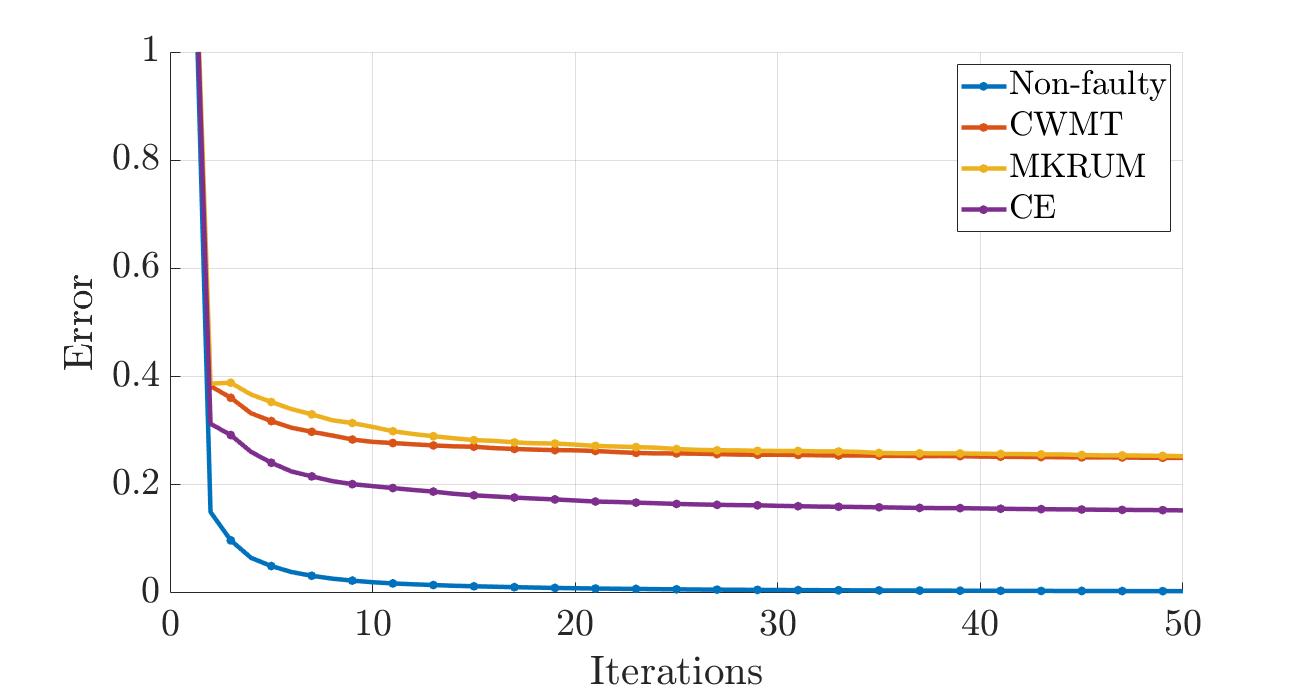}
        \label{fig:tiger}}
    \caption{The plots show the error $q^{\mathcal{H}}(\Bar{x}_{k}) -q^{\mathcal{H}}(x^{*})$ for $\mathcal{T}=3 $ for local iterations of local SGD with CE filter (Algorithm 1), CWMT, and Multi-KRUM for problem \eqref{eq:regression}. Going clockwise we vary the byzantine agents as $f = 2, 5, 8,10$.}
\end{figure*}

For the evaluation of the federated local SGD with CE filter, we consider a network with $N=50$ agents with a varying number of byzantine agents. We further compare the performance of the convergence of the algorithm with various other existing byzantine filters, namely multi-KRUM \cite{blanchard2017machine} and Coordinate-Wise Trimmed Median (CWTM) \cite{yin2018byzantine}, \cite{su2019finite}. Our experiment goals are the following:

\begin{itemize}
    \item Fix the number of byzantine agents compare the performace of the algorithm with differnt byzantine filters for $\mathcal{T} = 3$ local local and $50$ global communication rounds.
    \item For  $\mathcal{T} = 3$ local local and $50$ global communication rounds compare the performance of local SGD with CE filter with varying number of byzantine agents, $f = 2,5,8,10$.
\end{itemize}

We consider a scenario where we have $50$ agents trying to optimize sum of given local functions. We present results from Algorithm \ref{alg:cap} and local SGD with the aforementioned byzantine filters. Here at each local iteration any agent $i$ has access to an i.i.d sample of it's local gradients.
First, we consider a regression problem
\begin{align}\label{eq:regression}
    &\min_{x}\sum_{i=1}^{N}q^{i}(x)\triangleq \sum_{i=1}^{N}\Big(\|\mathbf{A}^{i}x - b^{i}\|^{2} + \sin^{2}(\|\mathbf{A}^{i}x - b^{i}\|)\Big).
\end{align}
where $\mathbf{A}_{i}$ and $b_{i}$ are the feature vector and labels respectively for $i^{\text{th}}$ agent.  This is an example of an invex but non-convex function satisfying the PL condition. Here each agent has its own estimate $x_{i}\in \mathbb{R}^{d}$ and also maintains the parameters $(\mathbf{A}_{i},b_{i})$ with $\mathbf{A}_{i} \in \mathbb{R}^{l\times d}$ and $b_{i} \in \mathbb{R}^{l}$. Further, each agent has an associated local function $q_{i}(x)$ where $x$ is the decision variable, $\mathbf{A}_{i}$ is the weight or importance of agent $i$'s objective, and $b_i$ is the target or reference value of user $i$'s objective. Here we further note that the since objectives satisfies the PL condition, the optimal solution $x^{*}$ will be unique for any set of honest agent $\mathcal{H}$ (see Remark 1 in \cite{gupta2023byzantine}). 

Second, we consider 

\begin{align}\label{eq:classification}
     \min_{x}\sum_{i=1}^{N}q_{i}(x)= \min_{x}\sum_{i=1}^{N}\Big(\frac{1}{1+\exp{(-\|\mathbf{A}_{i}x }-b_{i}\|)}\Big),
\end{align}
where the objectives are non-convex and do not satisfy the PL condition. Here we observe variation of the term $\frac{1}{|\mathcal{H}|}\sum_{i \in \mathcal{H}}\|\nabla q^{i}(\Bar{x}_{k})\|^2$.
This is a standard error term used for study analysis of a non-convex optimization problem. The simulation results are shown in Fig. 1 and 2 respectively. Fig.1 show the performances of CWTM, Multi-KRUM and CE filer with non-faulty Local GD as baseline for comparison as we vary the number of byzantine agents $f = 2,5,8,10$. We observe that Local SGD with CE filter out performs other byzantine filters. We further note that as the byzantine agent increase the convergence error increases. Fig 2 and 3 shows the performance of CE filter for local iterations $\mathcal{T} = 1,3$. We conclude that as we increase the local iterations, Algorithm \ref{alg:cap} converges faster which is consistent with previous work \cite{gupta2023byzantine}. The observed results are in line with our theoretical findings, which demonstrate that achieving only an approximate fault tolerance in stochastic scenarios is possible. 

\begin{figure*}[]\label{sin_2}
    \centering\subfigure{\includegraphics[width=0.45\linewidth, height=0.25\textwidth]{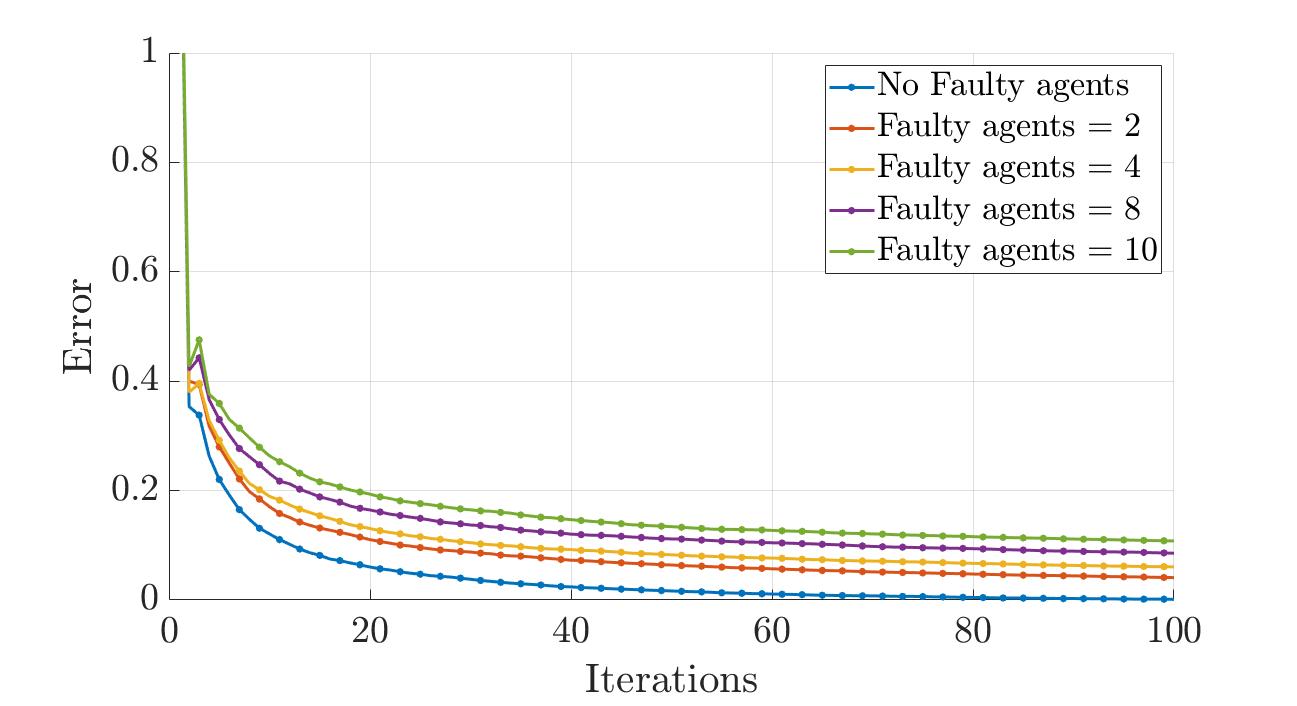}
        \label{fig:gull}}
     \subfigure{
        \includegraphics[width=0.45\linewidth, height=0.25\textwidth]{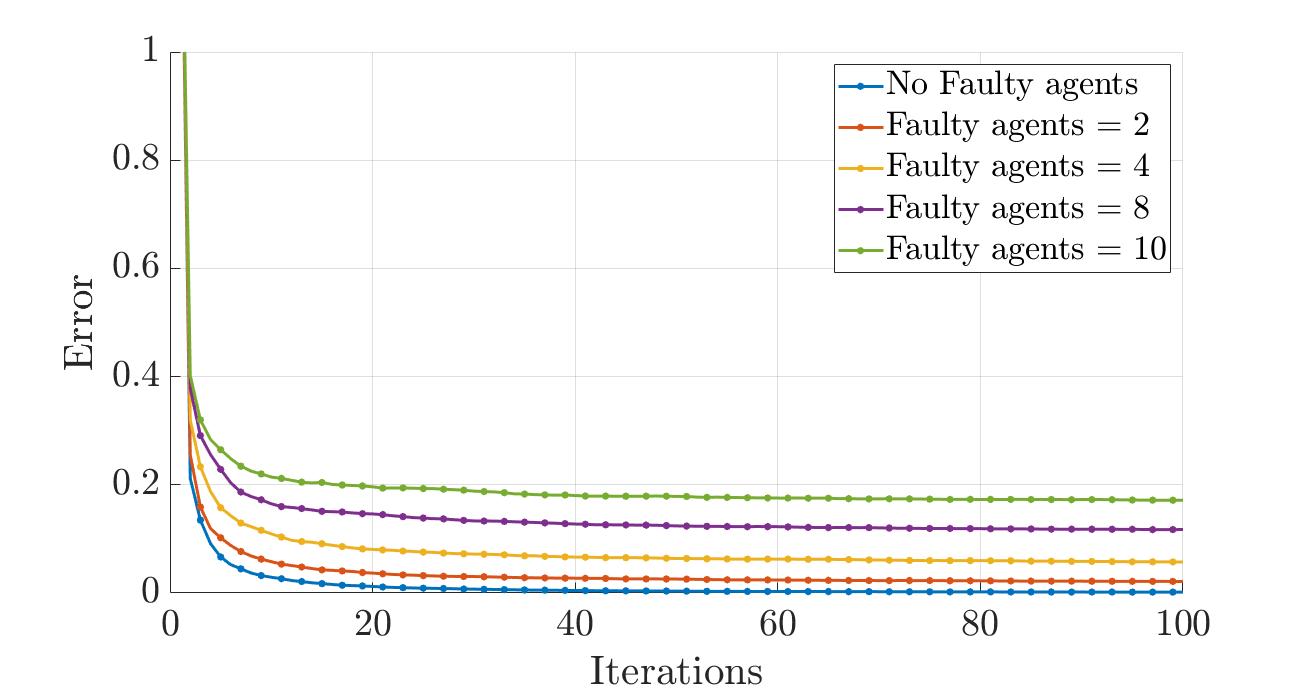}
        \label{fig:tiger}}
    \caption{The plots show the error $q^{\mathcal{H}}(\Bar{x}_{k}) -q^{\mathcal{H}}(x^{*})$ for problem \eqref{eq:regression} using Algorithm \ref{alg:cap}. Left and right figures show simulations with the local iterations $\mathcal{T}=1$ and $\mathcal{T}=3$, respectively.}
\end{figure*}

\begin{figure*}[]\label{sigmoid}
    \centering
     \subfigure{
        \includegraphics[width=0.45\linewidth, height=0.25\textwidth]{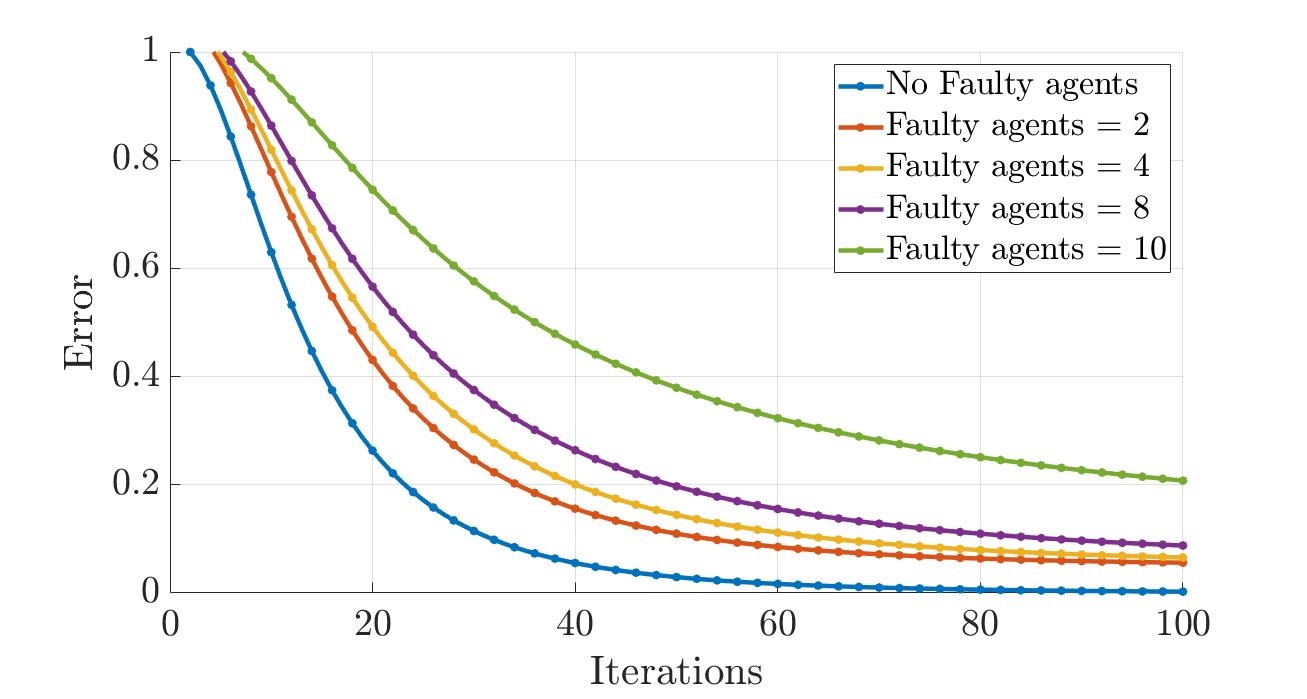}
        \label{fig:gull}}
     \subfigure{
        \includegraphics[width=0.45\linewidth, height=0.25\textwidth]{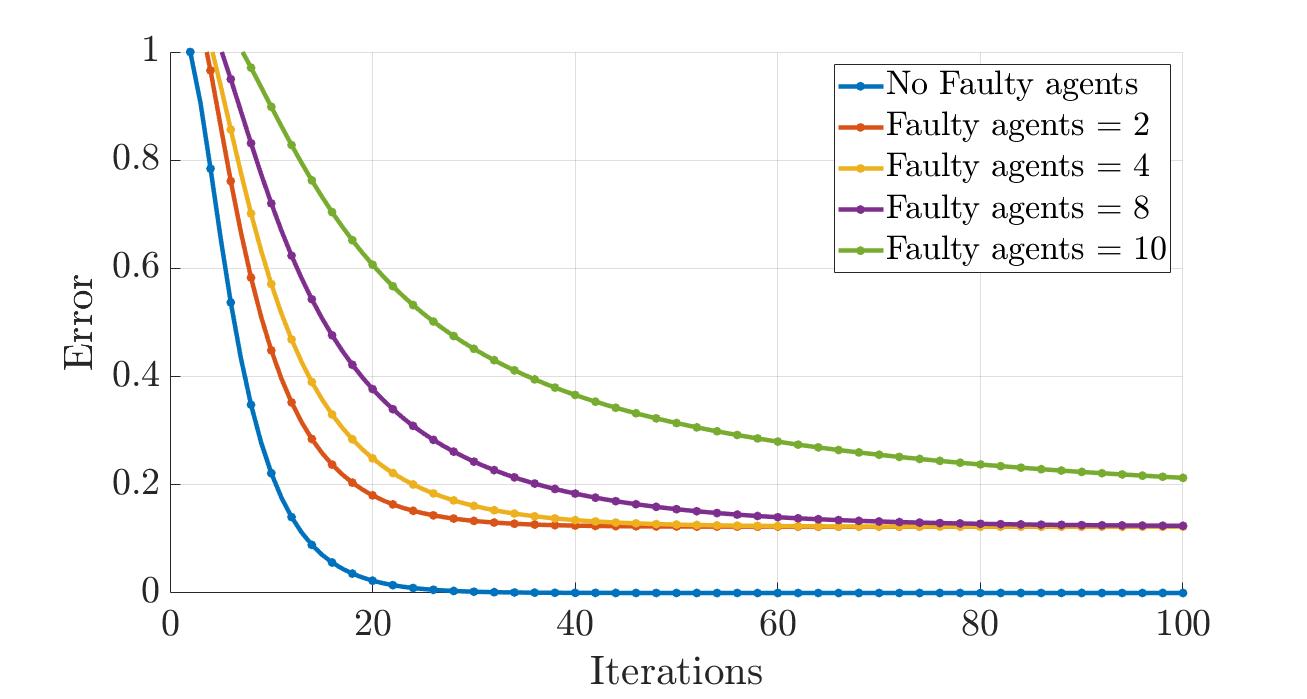}
        \label{fig:tiger}}
    \caption{The plots show the error $\frac{1}{|\mathcal{H}|}\sum_{i \in \mathcal{H}}\|\nabla q^{i}(\Bar{x}_{k})\|^2$ for local SGD with CE filter in solving problem \eqref{eq:classification}. Left and right figures show simulations with the local iterations $\mathcal{T}=1$ and $\mathcal{T}=3$, respectively.}
\end{figure*}

\section{Acknowledgement}
The work of Amit Dutta and Thinh T. Doan was partially supported by the Commonwealth Cyber Initiative, an investment in the advancement of cyber R\&D, innovation, and workforce development. For more information about CCI, visit \url{www.cyberinitiative.org}. 


\bibliographystyle{IEEEtran}
\bibliography{IEEEfull,ncs}

\begin{thebibliography}{10}
\providecommand{\url}[1]{#1}
\csname url@rmstyle\endcsname
\providecommand{\newblock}{\relax}
\providecommand{\bibinfo}[2]{#2}
\providecommand\BIBentrySTDinterwordspacing{\spaceskip=0pt\relax}
\providecommand\BIBentryALTinterwordstretchfactor{4}
\providecommand\BIBentryALTinterwordspacing{\spaceskip=\fontdimen2\font plus
\BIBentryALTinterwordstretchfactor\fontdimen3\font minus \fontdimen4\font\relax}
\providecommand\BIBforeignlanguage[2]{{%
\expandafter\ifx\csname l@#1\endcsname\relax
\typeout{** WARNING: IEEEtran.bst: No hyphenation pattern has been}%
\typeout{** loaded for the language `#1'. Using the pattern for}%
\typeout{** the default language instead.}%
\else
\language=\csname l@#1\endcsname
\fi
#2}}

\bibitem{mcmahan2017communication}
B.~McMahan, E.~Moore, D.~Ramage, S.~Hampson, and B.~A. y~Arcas, ``Communication-efficient learning of deep networks from decentralized data,'' in \emph{Artificial intelligence and statistics}.\hskip 1em plus 0.5em minus 0.4em\relax PMLR, 2017, pp. 1273--1282.

\bibitem{zhao2021fedpage}
H.~Zhao, Z.~Li, and P.~Richt{\'a}rik, ``Fedpage: A fast local stochastic gradient method for communication-efficient federated learning,'' \emph{arXiv preprint arXiv:2108.04755}, 2021.

\bibitem{woodworth2020local}
B.~Woodworth, K.~K. Patel, S.~Stich, Z.~Dai, B.~Bullins, B.~Mcmahan, O.~Shamir, and N.~Srebro, ``Is local sgd better than minibatch sgd?'' in \emph{International Conference on Machine Learning}.\hskip 1em plus 0.5em minus 0.4em\relax PMLR, 2020, pp. 10\,334--10\,343.

\bibitem{niknam2020federated}
S.~Niknam, H.~S. Dhillon, and J.~H. Reed, ``Federated learning for wireless communications: Motivation, opportunities, and challenges,'' \emph{IEEE Communications Magazine}, vol.~58, no.~6, pp. 46--51, 2020.

\bibitem{yang2020energy}
Z.~Yang, M.~Chen, W.~Saad, C.~S. Hong, and M.~Shikh-Bahaei, ``Energy efficient federated learning over wireless communication networks,'' \emph{IEEE Transactions on Wireless Communications}, vol.~20, no.~3, pp. 1935--1949, 2020.

\bibitem{shi2021challenges}
S.~Shi, Y.~Xiao, W.~Lou, C.~Wang, X.~Li, Y.~T. Hou, and J.~H. Reed, ``Challenges and new directions in securing spectrum access systems,'' \emph{IEEE Internet of Things Journal}, vol.~8, no.~8, pp. 6498--6518, 2021.

\bibitem{gupta2020fault}
N.~Gupta and N.~H. Vaidya, ``Fault-tolerance in distributed optimization: The case of redundancy,'' in \emph{Proceedings of the 39th Symposium on Principles of Distributed Computing}, 2020, pp. 365--374.

\bibitem{gupta2023byzantine}
N.~Gupta, T.~T. Doan, and N.~Vaidya, ``Byzantine fault-tolerance in federated local sgd under $2 f $-redundancy,'' \emph{IEEE Transactions on Control of Network Systems}, 2023.

\bibitem{blanchard2017machine}
P.~Blanchard, E.~M. El~Mhamdi, R.~Guerraoui, and J.~Stainer, ``Machine learning with adversaries: Byzantine tolerant gradient descent,'' \emph{Advances in neural information processing systems}, vol.~30, 2017.

\bibitem{su2019finite}
L.~Su and S.~Shahrampour, ``Finite-time guarantees for byzantine-resilient distributed state estimation with noisy measurements,'' \emph{IEEE Transactions on Automatic Control}, vol.~65, no.~9, pp. 3758--3771, 2019.

\bibitem{chen2017distributed}
Y.~Chen, L.~Su, and J.~Xu, ``Distributed statistical machine learning in adversarial settings: Byzantine gradient descent,'' \emph{Proceedings of the ACM on Measurement and Analysis of Computing Systems}, vol.~1, no.~2, pp. 1--25, 2017.

\bibitem{guerraoui2018hidden}
R.~Guerraoui, S.~Rouault, \emph{et~al.}, ``The hidden vulnerability of distributed learning in byzantium,'' in \emph{International Conference on Machine Learning}.\hskip 1em plus 0.5em minus 0.4em\relax PMLR, 2018, pp. 3521--3530.

\bibitem{li2019rsa}
L.~Li, W.~Xu, T.~Chen, G.~B. Giannakis, and Q.~Ling, ``Rsa: Byzantine-robust stochastic aggregation methods for distributed learning from heterogeneous datasets,'' in \emph{Proceedings of the AAAI Conference on Artificial Intelligence}, vol.~33, no.~01, 2019, pp. 1544--1551.

\bibitem{farhadkhani2022byzantine}
S.~Farhadkhani, R.~Guerraoui, N.~Gupta, R.~Pinot, and J.~Stephan, ``Byzantine machine learning made easy by resilient averaging of momentums,'' in \emph{International Conference on Machine Learning}.\hskip 1em plus 0.5em minus 0.4em\relax PMLR, 2022, pp. 6246--6283.

\bibitem{cao2019distributed}
X.~Cao and L.~Lai, ``Distributed gradient descent algorithm robust to an arbitrary number of byzantine attackers,'' \emph{IEEE Transactions on Signal Processing}, vol.~67, no.~22, pp. 5850--5864, 2019.

\bibitem{cao2020distributed}
------, ``Distributed approximate newton's method robust to byzantine attackers,'' \emph{IEEE Transactions on Signal Processing}, vol.~68, pp. 6011--6025, 2020.

\bibitem{yin2018byzantine}
D.~Yin, Y.~Chen, R.~Kannan, and P.~Bartlett, ``Byzantine-robust distributed learning: Towards optimal statistical rates,'' in \emph{International Conference on Machine Learning}.\hskip 1em plus 0.5em minus 0.4em\relax PMLR, 2018, pp. 5650--5659.

\end{thebibliography}

\section{appendix}\label{appendix}

\subsection{Proof of Theorem \ref{thm:Stochastic_Theorem2}}
\begin{proof}
Using the local SGD update for $i \in \Hcal$ and Assumptions \ref{eq:Assump_1}--\eqref{eq:Assump_3} we obtain the relations below. We skip their proofs due to space limitations.  
\begin{align}
\begin{aligned}
    &\mathbb{E}[\|x^{i}_{k,t+1}-x^{*}_{\Hcal}\|]{\leq}\frac{2}{\mu}\mathbb{E}[\|\nabla q^{\Hcal}(\Bar{x}_{k})\|]+\frac{2\sigma}{L},\\
    &\mathbb{E}[\|x^{i}_{k,t+1} - \Bar{x}_{k}\|]\leq\frac{2L\Tcal\alpha_{k}}{\mu}\mathbb{E}[\|\nabla q^{\Hcal}(\Bar{x}_{k})\|] + 3\sigma\Tcal\alpha_{k},\\
    &\mathbb{E}\big[\|x^{i}_{k,t+1}-x^{*}_{\Hcal}\|^{2}] \leq \frac{2}{\mu^{2}}\mathbb{E}[\|\nabla q^{\Hcal}(\Bar{x}_{k})\|^{2}] + \frac{2\sigma^{2}}{L^{2}}\cdot
\end{aligned}\label{eq:proof-thm1:ineq}
\end{align}
By \eqref{eq:Lipschitz}, we have
\begin{align}\label{eq:optimality_gap}
   \mathbb{E}[q^{\mathcal{H}}(\bar{x}_{k+1}) - q^{\mathcal{H}}(\bar{x}_{k})] &\leq  \mathbb{E}[\nabla q^{\mathcal{H}}(\bar{x}_{k})^{T}(\bar{x}_{k+1}-\bar{x}_{k})]\notag\\
   &\quad +\frac{L}{2}\mathbb{E}[\|\bar{x}_{k+1}-\bar{x}_{k}\|^{2}].
\end{align}
We next analyze each term on the right-hand side of \eqref{eq:optimality_gap}. Note that $|\Fcal_{k}|=|\Hcal|-N-f$ and $|\Bcal_{k}|=|\Hcal\backslash\Hcal_{k}|$. Thus, we have the following relation
\begin{align}
    \bar{x}_{k+1} 
    &=\frac{1}{|\Hcal|}\Big[\sum_{i \in \Hcal}x^{i}_{k,\Tcal}+\sum_{i \in \Bcal_{k}}x^{i}_{k,\Tcal}-\sum_{i \in \Hcal\backslash\Hcal_{k}}x^{i}_{k,\Tcal}\Big]\notag\\
    &\overset{\eqref{eq:Algo_T1}}{=}\Bar{x}_{k}-\frac{\alpha_{k}}{|\Hcal|}\sum_{i \in \Hcal}\sum_{t=0}^{\Tcal-1}\nabla Q^{i}(x^{i}_{k,t};\Delta^{i}_{k,t})+V_{x}\allowdisplaybreaks\notag\\
    &=  \bar{x}_{k} - \mathcal{T}\alpha_{k}\nabla q^{\Hcal}(\Bar{x}_{k})+ V_{x}\notag\\
    &\quad - \frac{\alpha_{k}}{|\mathcal{H}|}\sum_{i \in \mathcal{H}}\sum_{t=0}^{\mathcal{T}-1}(\nabla Q^{i}(x^{i}_{k,t};\Delta^{i}_{k,t}) -\nabla q^{i}(\Bar{x}_{k})).\notag
\end{align}
where 
\begin{align}
    V_{x} = \frac{1}{|\mathcal{H}|}\Big[ \sum_{i \in \mathcal{B}_{k}} (x^{i}_{k,\mathcal{T}}-\Bar{x}_{k}) - \sum_{i \in \mathcal{H} \backslash \mathcal{H}_{k}} (x^{i}_{k,\mathcal{T}}-\Bar{x}_{k})\Big].\label{eq:C_x}
\end{align}
We next consider the first term in \eqref{eq:optimality_gap} 
\begin{align}\label{eq:Term1}
    &\nabla q^{\mathcal{H}}(\bar{x}_{k})^{T}(\bar{x}_{k+1}-\bar{x}_{k})\notag\\
    &\leq-\frac{\alpha_{k}}{|\mathcal{H}|}\sum_{i \in \mathcal{H}}\sum_{t=0}^{\mathcal{T}-1}\nabla q^{\mathcal{H}}(\Bar{x}_{k})^{T}(\nabla Q^{i}(x^{i}_{k,t};\Delta^{i}_{k,t})- \nabla q^{i}(\Bar{x}_{k}))\notag\\
    &\quad -\mathcal{T}\alpha_{k}\|\nabla q^{\mathcal{H}}(\Bar{x}_{k})\|^{2} + \nabla q^{\mathcal{H}}(\Bar{x}_{k})^{T}V_{x}.
\end{align}
Analyzing the first term on the right-hand side of the above inequality we have,
\begin{align}
   & \frac{-\alpha_{k}}{|\mathcal{H}|}\sum_{i \in \mathcal{H}}\sum_{t=0}^{\mathcal{T}-1}\mathbb{E}[\nabla q^{\mathcal{H}}(\Bar{x}_{k})^{T}(\nabla Q^{i}(x^{i}_{k,t};\Delta^{i}_{k,t})
   \notag\\
   &-\nabla q^{i}(\Bar{x}_{k}))|\Pcal_{k,t}] \notag\\
   &=-\frac{\alpha_{k}}{|\Hcal|}\sum_{i \in \Hcal}\sum_{t=0}^{\Tcal-1}\nabla q^{\Hcal}(\Bar{x}_{k})^{T}(\nabla q^{i}(x^{i}_{k,t})-\nabla q^{i}(\Bar{x}_{k}))\notag\\
   &\leq \frac{L\alpha_{k}^{2}}{|\Hcal|}\sum_{i \in \Hcal}\sum_{t=0}^{\Tcal-1}\|\Bar{x}_{k}-x^{*}_{\Hcal}\|\|x^{i}_{k,t}-\Bar{x}_{k}\|\notag\\
   & \leq \frac{L^{3}\alpha^{2}_{k}\Tcal^{2}}{2\mu^{2}}\|\nabla q^{\Hcal}(\Bar{x}_{k})\|^{2} + \frac{L}{2\Tcal|\Hcal|}\sum_{i\in \Hcal}\sum_{t=0}^{\Tcal-1}\|x^{i}_{k,t}-\Bar{x}_{k}\|^{2},\notag
\end{align}
where the last inequality is due to Cauchy-Schwarz inequality $2xy \leq \eta x^{2}+y^{2}/\eta$ for any $\eta>0$ and $x,y \in \mathbb{R}$. Taking an expectation on both sides of the above inequality and using \eqref{eq:proof-thm1:ineq} we have 
\begin{align}
    &-\frac{\alpha_{k}}{|\mathcal{H}|}\sum_{i \in \mathcal{H}}\sum_{t=0}^{\mathcal{T}-1}\mathbb{E}[\nabla q^{\mathcal{H}}(\Bar{x}_{k})^{T}(\nabla Q^{i}(x^{i}_{k,t};\Delta^{i}_{k,t})- \nabla q^{i}(\Bar{x}_{k}))] \notag\\
    &\leq \frac{3L^{3}\alpha_{k}^{2}\Tcal^{2}}{2\mu^{2}}\mathbb{E}[\|\nabla q^{\Hcal}(\Bar{x}_{k})\|^{2}] + \frac{3L\Tcal^{2}\alpha^{2}_{k}\sigma^{2}}{2}.\notag
\end{align}
Next, we analyze $\|V_{x}\|^{2}$. Using \eqref{eq:C_x} and the fact that there exists an agent $j \in \mathcal{H}\backslash \Hcal_{k}$ such that $\|x^{i}_{k,t}-\Bar{x}_{k}\| \leq \|x^{j}_{k,t}-\Bar{x}_{k}\|$ for all agent $i \in \mathcal{B}_{k}$, we obtain
\begin{align}
    &\mathbb{E}\|V_{x}\|^{2}\notag\\
    &= \mathbb{E}\Big[\Big\|\frac{1}{|\mathcal{H}|}\big( \sum_{i \in \mathcal{B}_{k}} (x^{i}_{k,\mathcal{T}}-\Bar{x}_{k}) - \sum_{i \in \mathcal{H} \backslash \mathcal{H}_{k}} (x^{i}_{k,\mathcal{T}}-\Bar{x}_{k})\big)\Big\|^{2}\Big]\notag\\
    &\leq \frac{2|\mathcal{B}_{k}|^{2}}{|\mathcal{H}|^{2}}\mathbb{E}[\|x^{j}_{k,t}-\Bar{x}_{k}\|^{2}] +\frac{2|\mathcal{B}_{k}|}{|\mathcal{H}|^{2}}\sum_{i\in \mathcal{H} \backslash \mathcal{H}_{k}}\mathbb{E}[\|x^{i}_{k,t}-\Bar{x}_{k}\|^{2}], \notag\\
    & \overset{\eqref{eq:proof-thm1:ineq}}{\leq} \frac{8L^{2}\mathcal{T}^{2}\alpha_{k}^{2}}{\mu^{2}}\frac{|\mathcal{B}_{k}|^{2}}{|\mathcal{H}|^{2}}\mathbb{E}[\|\nabla q^{\mathcal{H}}(\Bar{x}_{k})\|^{2}] + \frac{12\mathcal{T}^{2}\alpha^{2}_{k}\sigma^{2}|\mathcal{B}_{k}|^{2}}{|\mathcal{H}|^{2}}.\notag
\end{align}
We now analyze the last term in the right-hand side of \eqref{eq:Term1}. Using the above result and the relation $<x,y>\eta\|x\|^{2}/2+\|y\|^{2}/2\eta \leq$ for any $\eta>0$ we have
\begin{align}
    &\mathbb{E}[\nabla q^{\mathcal{H}}(\Bar{x}_{k})^{T}V_{x}]\notag\\
    &\leq \frac{3L\Tcal\alpha_{k}|\Bcal_{k}|}{2\mu|\Hcal|}\mathbb{E}\big[\|\nabla q^{\Hcal}(\Bar{x}_{k})\|^{2}\big]  + \frac{\mu|\Hcal|}{6L\Tcal \alpha_{k}|\Bcal_{k}|}\mathbb{E}\big[\|V_{x}\|^{2}\big],\notag\\
    &\leq \frac{17L\Tcal\alpha_{k}|\Bcal_{k}|}{6\mu|\Hcal|}\mathbb{E}[\|\nabla q^{\Hcal}(\Bar{x}_{k})\|^{2}] + \frac{2\Tcal\mu\alpha_{k}\sigma^{2}|\Bcal_{k}|}{L|\Hcal|}\cdot\notag 
\end{align}
Thus we obtain from \eqref{eq:Term1} 
\begin{align}\label{eq:Term1_final}
    &\mathbb{E}[\nabla q^{\mathcal{H}}(\bar{x}_{k+1})^{T}(\bar{x}_{k+1}-\bar{x}_{k})] \notag\\
    &\leq \Big(-\alpha_{k}\mathcal{T} + \frac{17L\Tcal\alpha_{k}|\Bcal_{k}|}{6\mu|\Hcal|}+\frac{3L^{3}\alpha_{k}^{2}\Tcal^{2}}{2\mu^{2}}\Big)\mathbb{E}[\|\nabla q^{\mathcal{H}}(\Bar{x}_{k})\|^{2}]  \notag\\
    &\quad + \frac{3L\Tcal^{2}\alpha^{2}_{k}\sigma^{2}}{2}+ \frac{2\Tcal\mu\alpha_{k}\sigma^{2}|\Bcal_{k}|}{L|\Hcal|}\cdot
\end{align}
Next we analyze the second term on right hand side of \eqref{eq:optimality_gap}
\begin{align}
    &\|\Bar{x}_{k+1} - \Bar{x}_{k}\|^{2} \notag\\
    &= \Big\|\frac{\alpha_{k}}{|\Hcal|}\sum_{i \in \Hcal}\sum_{t=0}^{\Tcal-1}\nabla Q^{i}(x^{i}_{k,t};\Delta^{i}_{k,t})\Big\|^{2} +\|V_{x}\|^{2}\notag\\
    &\quad -\frac{2\alpha_{k}}{|\Hcal|}\sum_{i\in \Hcal}\sum_{t=0}^{\Tcal-1}V_{x}^{T}\nabla Q^{i}(x^{i}_{k,t};\Delta^{i}_{k,t}).\label{eq:Term2}
\end{align}
Taking the expectation of the first term on the right-hand side of \eqref{eq:Term2}
\begin{align}\label{eq:Term2_1}
    &\mathbb{E}\Big[\big\|\frac{\alpha_{k}}{|\Hcal|}\sum_{i \in \Hcal}\sum_{t=0}^{\Tcal-1}\nabla Q^{i}(x^{i}_{k,t};\Delta^{i}_{k,t})\big\|^{2}\Big]\notag\\
    &\leq \frac{\alpha^{2}_{k}\Tcal}{|\Hcal|}\sum_{i \in \Hcal}\sum_{t=0}^{\Tcal-1}\mathbb{E}\Big[\big\|\nabla Q^{i}(x^{i}_{k,t};\Delta^{i}_{k,t})\big\|^{2}\Big]\notag\\
    &\leq \frac{\alpha^{2}_{k}\Tcal}{|\Hcal|}\sum_{i \in \Hcal}\sum_{t=0}^{\Tcal-1}\mathbb{E}\Big[\big\|\nabla Q^{i}(x^{i}_{k,t};\Delta^{i}_{k,t})-\nabla q^{i}(x^{i}_{k,t})\big\|^{2}\Big]\notag\\
    &+\frac{\alpha^{2}_{k}\Tcal}{|\Hcal|}\sum_{i \in \Hcal}\sum_{t=0}^{\Tcal-1}\mathbb{E}\Big[\big\|\nabla q^{i}(x^{i}_{k,t})-\nabla q^{i}(x^{*}_{\Hcal})\big\|^{2}\Big]\notag\\
    &\leq \frac{L^{2}\Tcal\alpha_{k}^{2}}{|\Hcal|}\sum_{i \in \Hcal}\sum_{t=0}^{\Tcal-1}\mathbb{E}\Big[\big\|x^{i}_{k,t}-x^{*}_{\Hcal}\big\|^{2}\Big] + \Tcal^{2}\sigma^{2}\alpha^{2}_{k}\notag\\
    &\overset{\eqref{eq:proof-thm1:ineq}}{\leq}\frac{2L^{2}\Tcal^{2}\alpha^{2}_{k}}{\mu^{2}}\mathbb{E}\Big[\big\|\nabla q^{\Hcal}(\Bar{x}_{k})\big\|^{2}\Big] + 3\Tcal^{2}\sigma^{2}\alpha^{2}_{k}.
\end{align}
Now we analyze the last term in the right-hand side of \eqref{eq:Term2}. For this using the relation $<x,y>\eta\|x\|^{2}/2+\|y\|^{2}/2\eta \leq$ for any $\eta>0$ we have 
\begin{align}\label{eq:Term2_2}
    &\hspace{-0.2cm}\frac{-2}{|\Hcal|}\sum_{i\in \Hcal}\sum_{t=0}^{\Tcal-1}\mathbb{E}[\alpha_{k}V_{x}^{T}\nabla Q^{i}(x^{i}_{k,t};\Delta^{i}_{k,t})]\notag\\
    &\hspace{-0.2cm}\leq \frac{1}{|\Hcal|}\sum_{i\in \Hcal}\sum_{t=0}^{\Tcal-1}\Big(\frac{1}{\Tcal}\mathbb{E}[\|V_{x}\|^{2}] +\Tcal\alpha^{2}_{k}\mathbb{E}[\|\nabla Q^{i}(x^{i}_{k,t};\Delta^{i}_{k,t})\|^{2}]\Big)\notag\\
    &\hspace{-0.2cm}\leq\frac{\Tcal\alpha^{2}_{k}}{|\Hcal|}\sum_{i\in \Hcal}\sum_{t=0}^{\Tcal-1}\mathbb{E}[\|\nabla Q^{i}(x^{i}_{k,t};\Delta^{i}_{k,t})-\nabla q^{i}(x^{i}_{k,t})\|^{2}]\notag\\
    & +\frac{\Tcal\alpha^{2}_{k}}{|\Hcal|}\sum_{i\in \Hcal}\sum_{t=0}^{\Tcal-1}\mathbb{E}[\|\nabla q^{i}(x^{i}_{k,t})-\nabla q^{i}(x^{*}_{\Hcal})\|^{2}]+\mathbb{E}[\|V_{x}\|^{2}]\notag\\
    &\hspace{-0.2cm}\leq \mathbb{E}\|V_{x}\|^{2}+\frac{\Tcal L^{2}\alpha^{2}_{k}}{|\Hcal|}\sum_{i\in \Hcal}\sum_{t=0}^{\Tcal-1}\mathbb{E}[\|x^{i}_{k,t}-x^{*}_{\Hcal}\|^{2}] + \Tcal^{2}\sigma^{2}\alpha^{2}_{k}\notag\\
    &\hspace{-0.2cm}\overset{\eqref{eq:proof-thm1:ineq}}{\leq}\mathbb{E}\big[\|V_{x}\|^{2}\big]+\frac{2L^{2}\Tcal^{2}\alpha^{2}_{k}}{\mu^{2}}\mathbb{E}\big[\|\nabla q^{\Hcal}(\Bar{x}_{k})\|^{2}\big] + 3\Tcal^{2}\sigma^{2}\alpha^{2}_{k}.
\end{align}
Substituting \eqref{eq:Term2_1} and \eqref{eq:Term2_2} into \eqref{eq:Term2} we get obtain
\begin{align}\label{eq:Term2_final}
    &\frac{L}{2}\mathbb{E}[\|\Bar{x}_{k+1} - \Bar{x}_{k}\|^{2}]\notag\\
    &\leq 2\mathbb{E}\big[\|V_{x}\|^{2}\big] + \frac{4L^{2}\Tcal^{2}\alpha^{2}_{k}}{\mu^{2}}\mathbb{E}\big[\|\nabla q^{\Hcal}(\Bar{x}_{k})\|^{2}\big] + 6\Tcal^{2}\sigma^{2}\alpha^{2}_{k}\notag\\
    &\leq \Big(\frac{2L^{3}\Tcal^{2}\alpha^{2}_{k}}{\mu^{2}}+\frac{4L^{3}\mathcal{T}^{2}\alpha_{k}^{2}}{\mu^{2}}\frac{|\mathcal{B}_{k}|^{2}}{|\mathcal{H}|^{2}}\Big)\mathbb{E}[\|\nabla q^{\mathcal{H}}(\Bar{x}_{k})\|^{2}] \notag\\
    &\quad + 3L\Tcal^{2}\sigma^{2}\alpha^{2}_{k}+ \frac{6L\mathcal{T}^{2}\alpha^{2}_{k}\sigma^{2}|\mathcal{B}_{k}|^{2}}{|\mathcal{H}|^{2}}.
\end{align}
Finally, substituting \eqref{eq:Term1_final} and \eqref{eq:Term2_final} into \eqref{eq:optimality_gap} and using $|\Bcal_{k}|\leq f$ we get 
\begin{align}
    &\mathbb{E}[q^{\mathcal{H}}(\bar{x}_{k+1}) - q^{\mathcal{H}}(\bar{x}_{k})]\notag\\
    &\leq \Big(-\Big(1-\frac{17Lf}{6\mu|\Hcal|}\Big)\alpha_{k}\Tcal+ \frac{7L^{3}\Tcal^{2}\alpha^{2}_{k}}{2\mu^{2}}\notag\\
    &\quad +\frac{4L^{3}\Tcal^{2}\alpha_{k}^{2}f^{2}}{\mu^{2}|\Hcal|^{2}}\Big)\mathbb{E}\|\nabla q^{\mathcal{H}}(\Bar{x}_{k})\|^{2} + \frac{9L\Tcal^{2}\sigma^{2}\alpha^{2}_{k}}{2}\notag\\
    &\quad +\frac{2\Tcal^{2}\mu\alpha_{k}\sigma^{2}f}{L|\Hcal|} + \frac{6L\Tcal^{2}\sigma^{2}\alpha^{2}_{k}f^{2}}{|\Hcal|^{2}}\cdot\notag
\end{align}
Since $\frac{f}{|\mathcal{H}|} \leq \frac{\mu}{3L}$, we obtain \eqref{eq:assumption_1_Theorem2} 
\begin{align*}
    1-\frac{17Lf}{6|\Hcal|}>\frac{1}{18},
\end{align*}
using which we have
\begin{align}
    &\mathbb{E}[q^{\mathcal{H}}(\bar{x}_{k+1}) - q^{\mathcal{H}}(\bar{x}_{k})]\notag\\
    &\leq \Big(\frac{-\alpha_{k}\Tcal}{18}+ \Big(\frac{7}{2}+\frac{4f^{2}}{|\Hcal|^{2}}\Big)\frac{L^{2}\Tcal^{2}\alpha^{2}_{k}}{\mu^{2}}\Big)\mathbb{E}\big[\|\nabla q^{\mathcal{H}}(\Bar{x}_{k})\|^{2}\big] \notag\\
    &+ 5L\Tcal^{2}\sigma^{2}\alpha^{2}_{k} + \frac{2\Tcal^{2}\alpha_{k}\mu\sigma^{2}f}{L|\Hcal|}\cdot\notag
\end{align}
Next, using $\alpha_{k} = \alpha \leq \frac{\mu^{2}}{72L^{3}\mathcal{T}}$ from \eqref{eq:assumption_2_Theorem2}  we have,
\begin{align}
    &\mathbb{E}[q^{\mathcal{H}}(\bar{x}_{k+1}) - q^{\mathcal{H}}(\bar{x}_{k})] \notag\\
    &\leq -\frac{\alpha_{k}\mathcal{T}}{36}\mathbb{E}\big[\|\nabla q^{\mathcal{H}}(\Bar{x}_{k})\|^{2}\big] + 5L\Tcal^{2}\sigma^{2}\alpha^{2}_{k} + \frac{2\Tcal^{2}\alpha_{k}\mu\sigma^{2}f}{L|\Hcal|},\notag
\end{align}
which gives us 
\begin{align}
    &\mathbb{E}[q^{\mathcal{H}}(\bar{x}_{k+1}) - q^{\mathcal{H}}(x^{*}_{\mathcal{H}})] \notag\\
    &\leq \Big(1-\frac{\alpha_{k}\mu\mathcal{T}}{36}\Big)\mathbb{E}[q^{\mathcal{H}}(\bar{x}_{k}) - q^{\mathcal{H}}(x^{*}_{\mathcal{H}})]\notag\\
    &\quad + 5L\Tcal^{2}\sigma^{2}\alpha^{2}_{k} + \frac{2\Tcal^{2}\alpha_{k}\mu\sigma^{2}f}{L|\Hcal|}, \notag \\
    & \leq \Big(1-\frac{\alpha_{k}\mu\mathcal{T}}{36}\Big)^{k+1}\mathbb{E}[q^{\mathcal{H}}(\bar{x}_{0}) - q^{\mathcal{H}}(x^{*}_{\mathcal{H}})] \notag\\
    &\quad +\frac{180L\Tcal\alpha_{k}\sigma^{2}}{\mu}  + \frac{72\Tcal\sigma^{2}f}{\mu|\Hcal|}\cdot
\end{align}
This concludes our proof.
\end{proof}

\end{document}